\documentclass{article}

\usepackage{amsmath, amsthm, amscd, amsfonts, amssymb, graphicx, color}
\usepackage[bookmarksnumbered, colorlinks, plainpages]{hyperref}

\theoremstyle{definition}

\theoremstyle{remark}

\numberwithin{equation}{section}

\newcommand{\C}{\mathbb{C}}
\newcommand{\R}{\mathbb{R}}
\newcommand{\N}{\mathbb{N}}
\newcommand{\T}{\mathbb{T}}
\newcommand{\Z}{\mathbb{Z}}
\newcommand{\D}{\mathbb{D}}

\renewcommand{\ker}{\operatorname{ker}}
\newcommand{\ran}{\operatorname{ran}}

\newcommand{\dom}{\operatorname{dom}}
\newcommand{\norm}[1]{\parallel\!\!#1\!\!\parallel}
\renewcommand{\d}{\operatorname{d}}
\newcommand{\e}{\operatorname{e}}
\renewcommand{\i}{\operatorname{i}}

\newcounter{envcount}%

\newenvironment{Def}%
{\vspace{\bigskipamount}\refstepcounter{envcount}\textbf{(\theenvcount)\enspace Definition.}}%
  {\vspace{\bigskipamount}}

{\vspace{\bigskipamount}\enspace \textbf{Definition.}}%
  {\vspace{\bigskipamount}}

\newenvironment{The}%
{\vspace{\bigskipamount}\refstepcounter{envcount}\textbf{(\theenvcount)\enspace Theorem.}\itshape}%
  {\vspace{\bigskipamount}\upshape}

{\vspace{\bigskipamount}\refstepcounter{envcount}\textbf{(\theenvcount)\enspace Proposition.}\itshape}%
  {\vspace{\bigskipamount}\upshape}

\newenvironment{Cor}%
{\vspace{\bigskipamount}\refstepcounter{envcount}\textbf{(\theenvcount)\enspace Corollary.}\itshape}%
  {\vspace{\bigskipamount}\upshape}

{\vspace{\bigskipamount}\refstepcounter{envcount}\textbf{(\theenvcount)\enspace Lemma.}\itshape}%
  {\vspace{\bigskipamount}\upshape}

{\vspace{\bigskipamount}\refstepcounter{envcount}\textbf{(\theenvcount)\enspace}}%
  {\vspace{\bigskipamount}}

\theoremstyle{definition}
\swapnumbers

\begin{document}
\setcounter{page}{1}
\pagenumbering{arabic}

revised version august 21 for hrpub
\begin{center}
{\Large Unbounded Toeplitz Operators with Rational Symbols}\\  

\vspace{1cm}
Domenico P.L. Castrigiano\\
Technischen Universit\"at M\"unchen, Fakult\"at f\"ur Mathematik, M\"unchen, Germany\\

\smallskip

{\it E-mail address}: {\tt
castrig\,\textrm{@}\,ma.tum.de
}
\end{center}

\begin{quote} Unbounded (and bounded) Toeplitz operators (TO) with rational symbols are analysed in detail showing that they are densely defined closed and have finite dimensional kernels and deficiency spaces. The latter spaces as well as the domains, ranges, spectral and Fredholm points are  determined.  In particular, in the symmetric case, i.e., for a real rational symbol the deficiency spaces and indices are explicitly available. --- The concluding section gives a brief overview on the research on unbounded TO in order to locate the present contribution. Regarding properties of unbounded TO in general,  it furnishes some new results recalling the close relationship to Wiener-Hopf operators and, in case of semiboundedness, to singular operators of Hilbert transformation type. Specific symbols considered in the literature admit further analysis. Some conclusions are drawn for semibounded integrable and real square-integrable symbols. There is an approach to semibounded TO, which starts from  closable semibounded  forms related to a Toeplitz matrix. The Friedrichs extension of the TO associated with such a form is studied. Finally, analytic TO and Toeplitz-like operators are briefly examined, which in general differ from the TO treated here.
\\

{\it Mathematics Subject Classification MSC2020}: 47B35, 47A53, 47B25 
\end{quote}

\section{Introduction} Toeplitz matrices and bounded Toeplitz operators are thoroughly studied  for a  long time \cite{BS90}, and there are also important early results on unbounded Toeplitz operators as  \cite{R62} and \cite{R65} regarding their spectral representation in case of real semibounded integrable symbols.  More recently there is an increasing interest in unbounded Toeplitz and Toeplitz-like operators. 
So see e.g.\,\cite{H88},\,\cite{S08} for analytic Toeplitz operators and their domains and in particular for analytic Toeplitz operators with rational symbols, \cite{Y18} for closable quadratic forms for semibounded Toeplitz operators, and \cite{GHJR18} for a kind of  Toeplitz operators with rational symbols. For a brief overview on unbounded Toeplitz operators and some further details see sec.\,\ref{RTO}.
\\
\hspace*{6mm} 
Future research will  aim to extend classical result on bounded Toeplitz operators to the unbounded case. Here the class of Toeplitz operators with rational symbol is particularly useful as being relatively easy to handle.
\\
\hspace*{6mm} 
In general unbounded Toeplitz operators can be approached like unbounded Wiener-Hopf operators in \cite{C20} making use of  analogous methods and tools, in particular coming from complex function theory. Moreover, in the Hilbert space case, there is a well-known isomorphism connecting  Wiener-Hopf  and Toeplitz operators. The present article uses this connection and is essentially an application of the results in  \cite{C20} on Wiener-Hopf operators with rational symbols. Hence the proofs will be brief. 
\\
\hspace*{6mm} 
Let the torus $\T$ be endowed with the normalized Lebesgue measure $\lambda$. The Hardy space $H^2(\T)$ is the subspace of $L^2(\T)$ with orthonormal basis  
$e_n(z):=z^n$, $n\in\N_0$. Let $P_{H^2(\T)}: L^2(\T)\to H^2(\T)$ denote the orthogonal projection. Its adjoint $P^*_{H^2(\T)}$ is the identical injection of  $H^2(\T)$ into $L^2(\T)$.
For measurable  $\omega: \T\to \C$, $M(\omega)$ denotes the multiplication operator by $\omega$ in $L^2(\T)$.  

\begin{Def}\label{DTO}
The Toeplitz  operator (TO) $T_\omega$ in $H^2(\T)$ with symbol $\omega$ is defined by $$T_\omega :=P_{H^2(\T)}M(\omega)\,P^*_{H^2(\T)}$$
 So $T_\omega$ is the  trace (compression) of the multiplication operator $M(\omega)$  on the Hardy space just as the well-studied TO for  a bounded symbol. As mentioned some  properties of unbounded TO for general symbols are listed in  sec.\,\ref{RTO}.    
\end{Def}

\section{Rational Symbols} 
TO for rational symbols $\omega=\frac{R}{S}\big|_{\T}$ with complex polynomials $R\ne 0$, $S\ne 0$ permit the following more detailed analysis.  Important is the case of poles on the unit circle, i.e., of zeros of $S$ on $\T$, due to which $T_\omega$ is unbounded.   The analysis concerns the domains, ranges, kernels, and deficiency spaces and determines the Fredholm points  and  various parts of the spectrum. Recall that a densely defined closed operator between Banach spaces with finite dimensional kernel  and cokernel is called a Fredholm operator if its range is closed  (cf.\,\cite{S67}). 
\\
\hspace*{6mm}
 Mostly we will omit $|_\T$, which indicates the restriction on $\T$. By definition a polynomial with a negative degree is the null function.
 \\
\hspace*{6mm} 
Let the polynomials $R$,\,$S$ have no common zeros. Write $S=s\,S_{in}\,S_{ex}$, where the zeros of the polynomials $s$, $S_{in}$, and $S_{ex}$   are  the zeros of $S$ with modulus $=1$, $<1$, and $>1$, respectively. Put also $S_{\overline{ex}}:=s\,S_{ex}$, $S_{\overline{in}}:=s\,S_{in}$. Finally, let $\breve{S}$ denote the polynomial, whose    zeros are $1/\overline{b}$  with the same multiplicity as $b$, where $b\ne 0$ is a zero of $S$. So the zeros of $\breve{S}_{\,ex}$ are $1/\overline{b}$ for $b$ zero of $S$ with $0<|b|<1$.  Analogous notations concern $R$.

\begin{The}\label{RSKCK} Let $\omega=\frac{R}{S}$ be a rational function with polynomials $R$ and $S$ without common zeros.  Then $T_\omega$ is densely defined and closed and
\begin{itemize}
\item[\emph{(a)}] $\dom  T_\omega=s\, H^2(\T)$\; and \;$\ran T_\omega=P_{H^2(\T)}\big(\frac{R}{S_{in}S_{ex}}H^2(\T)\big)$
\item[\emph{(b)}] $\ker T_\omega=\Big\{S_{\overline{ex}}\,t/R_{ex}: \;t\textrm{ polynomial with }\deg t< \deg S_{in}-\deg R_{\overline{in}} \,\Big\}$
\item[\emph{(c)}] $\big(\ran T_\omega\big)^\perp=\Big\{\breve{S}_{\,ex}\,t/\breve{R}_{\,ex}:  \;t\textrm{ polynomial with }\deg t< \deg R_{in}-\deg S_{in}\Big\}$
\item[\emph{(d)}] the following statements are equivalent:
\begin{itemize}
\item[\emph{(1)}] $T_\omega$ is a Fredholm Operator
\item[\emph{(2)}] $\ran T_\omega$  closed 
\item[\emph{(3)}] $R$ has no zeros on $\T$, i.e., $\omega$ has no zeros     
\end{itemize}
\end{itemize}
\end{The}
 {\it Proof.} The proof consists in  transferring  the assertions to  Hilbert space isomorphic Wiener-Hopf operators and applying \cite[sec.\,3]{C20}. Indeed, as recalled in \cite[(2.3)]{C20}, there is  the well-known relation
\begin{equation}\label{THS}
T_\omega=\Gamma^{-1} M_+(\omega\circ C)\, \Gamma
\end{equation}
where $\Gamma:H^2(\T) \to H^2(\R)$, $(\Gamma u)(x):=\frac{\i}{\sqrt{\pi}(x+\i)}u\big(C(x)\big)$ is the Hilbert space isomorphism based on the Cayley transformation  $C(z):=\frac{z-\i}{z+\i}$, and $M_+(\kappa)=P_{H^2(\R)}M(\kappa)P^*_{H^2(\R)}$ is the Wiener-Hopf operator in $H^2(\R)$ for the measurable symbol $\kappa:\R\to \C$ with $P_{H^2(\R)}: L^2(\R)\to H^2(\R)$  the orthogonal projection and its adjoint $P^*_{H^2(\R)}$  the identical injection of  $H^2(\R)$ into $L^2(\R)$. We treat the case that $1$ is not a zero of $R\,S$. The general case follows analogously using for convenience  the transformation $\e^{\i\alpha}C$ in place of $C$ with $\e^{\i\alpha}$ not a zero of $R\,S$. \\
\hspace*{6mm}
For $\kappa:=\omega\circ C$  with $\omega=R/S$ one  readily finds $\kappa=P/Q$ with  polynomials $P,Q$ without common zeros given by 
\begin{itemize}
\item[(i)] $P(z):=a\,(z+\i)^{\max\{0,n-m\}}\prod_{j=1}^m\big(z-C^{-1}(a_j)\big)$ 
\item[(ii)] $Q(z):=b\,(z+\i)^{\max\{0,m-n\}}\prod_{k=1}^n\big(z-C^{-1}(b_k)\big)$ 
\end{itemize}
with constant factors
$a\ne 0$, 
$b\ne 0$.
Here $(a_j), (b_k)$ are the zeros of $R$ and $S$, respectively.\\
\hspace*{6mm}
 Exemplarily we show (c). First note $(\ran T_\omega)^\perp= \Gamma^{-1}(\ran M_+(\kappa))^\perp$. By \cite[3.2(c)]{C20}, $\big(\ran M_+(\kappa)\big)^\perp=\big\{(\tilde{Q}_</\tilde{P}_<)\,r:  \;r\textrm{ polynomial with }\deg r< \deg P_>-\deg Q_>\big\}$. As to the notation, given a polynomial $Y$, then the zeros of the polynomial  $Y_>$ ($Y_<$)  are the zeros of $Y$ in the upper  (lower) half-plane, and the coefficients of $\tilde{Y}$  are the complex conjugates ones of $Y$. One has to evaluate  the right-hand side of
 $$\Gamma^{-1}\big((\tilde{Q}_</\tilde{P}_<)\,r\big)(z)=2\sqrt{\pi}(1-z)^{-1}\big((\tilde{Q}_</\tilde{P}_<)\,r\big)\big(C^{-1}(z)\big)$$
 Now  $\tilde{P}_<\big(C^{-1}(z)\big)=\overline{a}\prod_{\{j:|a_j|<1\}} \big(C^{-1}(z)+C^{-1}(\overline{a}_j)\big)$, since $\overline{C^{-1}(a_j)}=-C^{-1}(\overline{a}_j)$. Note $C^{-1}(z)+C^{-1}(\overline{a}_j)=\frac{2\i}{1-\overline{a}_j}(1-z)^{-1}(1-\overline{a}_jz)$. Hence $\tilde{P}_<\big(C^{-1}(z)\big)=A\,(1-z)^{-\deg R_{in}}\breve{R}_{ex}$
  with $A\ne 0$. Similarly $\tilde{Q}_<\big(C^{-1}(z)\big)=B\,(1-z)^{-\deg S_{in}}\breve{S}_{ex}$
  with $B\ne 0$. Furthermore,  for $ l\le\deg P_>-\deg Q_>-1=\deg R_{in}-\deg S_{in}-1$ and $\deg r=l$ one finds $r\big(C^{-1}(z)\big)=(1-z)^{-l}t(z)$. Here $t$ is a polynomial with $t(1)\ne 0$  and $\deg t=l-\mu$, where $0\le\mu\le l$ is the multiplicity of the zero $-\i$ of $r$. Putting these results together proves (c).
\\  
\hspace*{6mm}
As to the proof of (d), note that $\deg P=\deg Q$ and that $P$ has no real zero if and only if $R$ has no zero on $\T$. 
 Apply  \cite[3.2(d)]{C20}.  \qed

Clearly, $T_\omega$ is bounded if and only if $S$ has no zeros on $\T$. For bounded TO recall the book \cite{BS90}. From  (\ref{RSKCK}) one immediately obtains

\begin{Cor}\label{SWHFP} Let $\lambda\in\C$ and put $R^\lambda:=R+\lambda S$. Then referring to $T_\omega$, $\lambda$ is 
\begin{itemize}
\item[\emph{(a)}] a Fredholm point  (i.e. $T_\omega-\lambda I$ is a Fredholm operator) iff $\lambda\not\in \omega(\T)$; if $\lambda$ is a Fredholm point, then $\dim\ker(T_\omega-\lambda I)=\max\{0,\deg S_{in}-\deg R^\lambda_{in}\}$,
 $\dim\ran(T_\omega-\lambda I)^\perp=\max\{0,\deg R^\lambda_{in}-\deg S_{in}\}$, and $\operatorname{ind}(T_\omega-\lambda I)=\deg S_{in}-\deg R^\lambda_{in}$
\item[\emph{(b)}]  a regular value (i.e. $T_\omega-\lambda I$ is continuously invertible) iff $\lambda\not\in \omega(\T)$ and $\deg S_{in}\le \deg R^\lambda_{in}$
\item[\emph{(c)}]   in the resolvent set iff $\lambda\not\in \omega(\T)$ and $\deg S_{in}= \deg R^\lambda_{in}$
\item[\emph{(d)}]  a spectral value iff $\lambda\in  \omega(\T)$ or  $\deg S_{in}\ne \deg R^\lambda_{in}$
\item[\emph{(e)}]  in the point spectrum iff \,$\deg R^\lambda_{\overline{in}} <\deg S_{in}$
\item[\emph{(f)}]  in the continuous spectrum (i.e.  $T_\omega-\lambda I$ is injective with dense not closed range)  iff $\lambda\in  \omega(\T)$ and $\deg R^\lambda_{in} \le\deg S_{in}\le \deg R^\lambda_{\overline{in}}$
\item[\emph{(g)}]  in the residual spectrum (i.e.  $T_\omega-\lambda I$ is injective with not dense range) iff $\deg S_{in} <\deg R^\lambda_{in}$
\end{itemize}
\end{Cor}

Note that $\lambda\mapsto \deg R^\lambda_{in}$ is locally non-decreasing  due to the continuity of the roots of a polynomial on its coefficients \cite{CC89}. On  $\C\setminus \omega(\T)$ it is even locally constant, since there
 $R^\lambda_{in}=R^\lambda_{\overline{in}}$. Hence, besides $\operatorname{ind}(T_\omega-\lambda I)$, also
$\dim\ker(T_\omega-\lambda I)$ and
 $\dim\ran(T_\omega-\lambda I)^\perp$ are constant on the components of $\C\setminus \omega(\T)$.\\

A densely defined TO is symmetric if and only its symbol is  real a.e. (by \cite[(2.3),\,2.14]{C20}).

\begin{Cor}\label{STO} Let $\omega$ be a real-valued rational function, i.e., $\omega(\T)\subset \R$ and $\omega=\frac{R}{S}$ with polynomials  $R,S$ without common zeros. Put
$$R-\i S=p\;q$$
where the zeros of the polynomials $p$ and $q$ have  modulus $<1$ and $>1$, respectively, and put $l:=\deg(R+\i S)-\deg\breve{p}-\deg\breve{q}$. \\
\hspace*{6mm}
Then $T_\omega$  is densely defined  closed symmetric. For $\dom T_\omega$ and $\ran T_\omega$ recall \emph{(\ref{RSKCK})(a)}. The deficiency  indices are
\begin{equation}
  n_+:=\dim \ran\big(T_\omega - \i I\big)^\perp= \deg p-\deg S_{in},\quad n_-:=\dim \ran\big(T_\omega + \i I\big)^\perp=l+ \deg q-\deg S_{in}
  \end{equation}
and the deficiency spaces  are
\begin{equation}
\ran\Big(T_\omega- \i I\Big)^\perp=\Big\{\frac{\breve{S}_{\,ex}}{\breve{p}}t: t\text{ polynomial with }\deg t<\deg p-\deg S_{in}\Big\}
 \end{equation}
 \begin{equation}
\ran\Big(T_\omega+ \i I\Big)^\perp=\Big\{\frac{\breve{S}_{\,ex}}{q}t: t\text{ polynomial with }\deg t< l+\deg q -\deg S_{in}\Big\}
 \end{equation}
\end{Cor}\\
{\it Proof.} Check that  $R-\i S$ has no zeros on $\T$.  Note $T_\omega - \i I=T_{\omega_+}$ with $\omega_+:=p\,q/S$. Since $(p\,q)\,\breve{}_{\,ex}=\breve{p}_{\,ex}$, the claim regarding   $\ran (T_\omega-\i I)^\perp$ and $n_+$ follows  from (\ref{RSKCK})(c) and \cite[3.4]{C20}.
\\
\hspace*{6mm}
 The result about the other deficiency space follows analogously once having shown that $R+\i S =cz^l\breve{p}\,\breve{q}$ for  
some $c\ne 0$. Indeed, since $P$ and $Q$ in (\ref{RSKCK}) are real, one has $P-\i Q=Q_+Q_-$ with the zeros of $Q_+$ and $Q_-$ in the upper and lower half-plane, respectively. Then $P+\i Q=\tilde{Q}_+\tilde{Q}_-$ replacing the coefficients of $Q_\pm$ by the complex conjugate ones. The claim follows by straightforward computation  using $C(\overline{a})\overline{C(a)}=1$ for $a\ne\pm\i$.
 \qed \\
 
 For an example regarding (\ref{STO}) with $l=1$ and  $\deg S_{in}=0$, consider $\omega=\i (z-a)/(z+a)$ for $a\in\T$. Check that $\omega$ on $\T$ is real, and  $R-\i S=-2\i a$,  $R+\i S=2\i z$. 
 It follows $n_+=0$, $n_-=1$, and $\ran (T_\omega+ \i I )^\perp=\C e_0$. So $T_\omega$ is unbounded densely defined closed maximal symmetric not self-adjoint. --- Similarly, an example for a bounded self-adjoint $T_\omega$ in (\ref{STO}) is $\omega=(z^2+1)/z$.
 \\
\hspace*{6mm}
Generally, if  $\omega$ is holomorphic  on the disc $\mathbb{D}$,
  then $T_\omega$ in (\ref{STO}) is  self-adjoint only if $\omega$ is constant. (This follows easily from $n_+=n_-=0$,  $\deg S_{in}=0$.) --- Finally there is no unbounded self-adjoint   $T_\omega$ in (\ref{STO}). (Indeed, assume the contrary. Assume $\deg (R-\i S)=\max\{\deg R,\deg S\}$. Otherwise consider $2\omega$ in place of $\omega$. Now, $\deg S_{in}=\deg S_{ex}$ (since $\omega$ is real), there is a zero of $S$ on $\T$ (since $T_\omega$ is unbounded), and $k:=\deg p=\deg S_{in}$, $l+\deg q=\deg S_{in}$ (since $n_+=n_-=0$), whence the contradiction $2k-l=\deg p+\deg q=\deg (R-\i S)\ge \deg S\ge 2k+1$.)

\section{Remarks}\label{RTO}

We recall some earlier results on unbounded TO in order to locate the present contribution.  The paragraphs (e),\,(f)  below concern in particular rational symbols. 

\subsubsection*{(a)  Some general properties of  unbounded TO}

The Toeplitz operator $T_\omega$ for measurable  symbol $\omega: \T\to \C$ is Hilbert space isomorphic to the Wiener-Hopf operator $M_+(\kappa)$ in $H^2(\R)$ with symbol $\kappa = \omega\circ C$, which is the trace on the Hardy space $H^2(\R)$  of the multiplication operator by $\kappa$ in $L^2(\R)$. See  \cite[(2.3)]{C20} recalled in  (\ref{THS}). This isomorphism permits to transfer results on  $M_+(\kappa)$ in \cite[sec.\,2]{C20} to $T_\omega$. Note also   that in \cite{CMP21} one finds a detailed study of the kernel of $M_+(\kappa)$. Thus one has the following properties. \\
\hspace*{6mm}
 The domain  $\dom T_\omega=\{u\in H^2(\T): \omega u\in L^2(\T)\}$ is either trivial or dense. It is dense if and only if $\ln(1+|\omega|)$, or equivalently $\ln^+|\omega|$, is integrable, and it equals $H^2(\T)$ if and only if $\omega$ is essentially bounded or, equivalently, $T_\omega$ is bounded. \\
\hspace*{6mm} 
For all $\alpha\in H^\infty(\T)$ one has \,(1) \;$T_\omega T_\alpha=T_{\alpha \omega}$ and 
 \,(2) \;$T_\alpha^* T_\omega\subset T_{\overline{\alpha} \omega}$.   Hence the shift invariance $T_\omega=T_z^*T_\omega T_z$ holds. \\
 \hspace*{6mm}
 For $b\ge 0$, $T_{(b)}:=\Gamma^{-1}M_+(\e^{\i bx})\,\Gamma$ is the TO with symbol $z\to\exp(b\frac{z+1}{z-1})$ in $H^\infty(\T)$. 
  One has  the   invariance $T_\omega=T_{(b)}^*T_\omega T_{(b)}$,
which for bounded TO is even characteristic. Up to some technical assumptions this holds true also for unbounded TO (cf.\,\cite[2.17]{C20}).\\
\hspace*{6mm}
If $\omega$ is not almost constant and if $\lambda\in\C$ is an eigenvalue of $T_\omega$ then $\overline{\lambda}$ is not an eigenvalue of $T_\omega$. If  $\omega$ is almost real not almost constant, then $T_\omega$   has no eigenvalues and in particular $T_\omega$  is injective. \\
\hspace*{6mm}
 If $T_\omega$ is densely defined,  then $T_\omega$ is symmetric if and only if $\omega$ is almost real. \\
 \hspace*{6mm}
  Let $T_\omega$ be densely defined symmetric. Then $T_\omega$ is  bounded below if and only if $\omega$ is real  essentially bounded below, and the maximal lower bound  of  $T_\omega$ equals the maximal essential lower bound of $\omega$. If $T_\omega$ is bounded below and not bounded, then the numerical range $\{\langle u, T_\omega u\rangle: u\in\dom T_\omega,\;||u||=1\}$ equals $]\alpha,\infty[$ with $\alpha$ the maximal lower bound. \\
  \hspace*{6mm}
 Finally, every densely defined semibounded TO (i.e., $\omega$ real semibounded and  $\ln(1+|\omega|)$ integrable) has a  canonical representation by a product of a closable operator and its adjoint. Then the Friedrichs extension of the TO is obtained replacing the operator by its closure.  By this factorization one has a rather explicit expression of the former. Moreover, the factorization determines  a  Hilbert space isomorphism relating  the semibounded TO  to a singular integral operator of Hilbert transformation type. This is shown in  \cite[sec.\,5]{C20} for the corresponding Wiener-Hopf operator and 
illustrated by the example of Lalescu's operator  \cite[sec.\,5.2]{C20}. As a further example, the spectral representation of the TO with symbol $1_A$, where $A\subset  \T$ is an arc, is explicitly related to that of the finite Hilbert transformation \cite[sec.\,3.2.3, 3.3.2]{C19}.

\subsubsection*{(b) Semibounded integrable symbols} 

For integrable symbol $\omega$ the TO $T_\omega$ is densely defined  and, if in addition $\omega$ is real bounded below,  $T_\omega$ is symmetric bounded below. Note however that $e_n\in\dom T_\omega$ only if $\omega$ is even square-integrable. Note also that an unbounded rational symbol is not integrable.  Writing $\omega -\lambda>0$ for $\lambda\in\R$ as the absolute square of an outer function (which corresponds to the above mentioned  factorization) and   analyzing the Toeplitz matrix $\big(\int_\T\overline{e}_n\,\omega\, e_m\d \lambda\big)_{nm}$, in \cite{R62},\,\cite{R65} the spectral representation of the Friedrichs extension of $T_\omega$ is deduced. If $\omega$ is not almost constant, than the latter is absolutely continuous.\\
\hspace*{6mm}
 One recalls that the Friedrichs extension $\tilde{T}_\omega$ of $T_\omega$ is given by $T^*_\omega|_{\dom T^*_\omega\,\cap\,\mathcal{D}_\omega}$. Here $\mathcal{D}_\omega:=\{u\in H^2(\T):\exists\; u_n\in\dom T_\omega\textrm{ with }||u_n-u_m||_\omega\to 0,\;||u_n-u||\to 0\}$ is the completion of $\dom T_\omega$ with respect to the norm 
$\norm{\cdot}_\omega$ defined by $\norm{u}^2_\omega=\langle u,T_\omega u\rangle+(1-\gamma)\norm{u}^2$.

\subsubsection*{(c) Closable semibounded  form related to a Toeplitz matrix} 
Let $(\omega_k)_{k\in\Z}$  be a complex sequence with $\overline{\omega}_k=\omega_{-k}$ such that the sequilinear  form $s$, determined by the Toeplitz matrix $(\omega_{n-m})_{n,m\ge0}$, is bounded below (i.e., $s(\xi,\xi)=\sum_{n,m\ge0}  \overline{\xi}_n\omega_{n-m}\xi_m\ge \gamma \sum_n|\xi_n|^2$ for some $\gamma\in\R$ and all  $\xi\in\mathcal{D}_0$, the dense subspace of $l^2(\N_0)$ of all finite sequences). 
Then by \cite[2.3]{Y18}, the form $s$ is closable if and only if there is a real  integrable function $\omega$ on $\T$ with lower bound $\gamma$ so that $\omega_k=\int_{\T} z^{-k}\omega(z)\d\lambda(z)$, $k\in\Z$. Here $\omega$ is uniquely determined. \\
\hspace*{6mm}
Let us identify $l^2(\N_0)$ with $H^2(\T)$ by the usual Hilbert space isomorphism $\xi\to \sum_{n\ge0}\xi_ne_n$. Consider the case that $s$ is closable. Then by \cite[3.3]{Y18}, the closure of $s$,  still denoted by $s$, is given by $s(u,v)=\int_\T\overline{u}\,v\,\omega \d \lambda$ for all $u,v\in H^2(\T)\cap L^2(\T,\Lambda)$. Here  $\Lambda:=\rho\lambda$, where the density $\rho:=\omega+1-\gamma$ is $\ge 1$ and integrable.
\\
\hspace*{6mm}
 $\mathcal{D}_\Lambda:=H^2(\T)\cap L^2(\T,\Lambda)$ is dense in $H^2(\T)$ (since  $\mathcal{D}_0\subset\mathcal{D}_\Lambda$) and closed in $L^2(\T,\Lambda)$ (since $||\cdot||_\Lambda\ge||\cdot||$). By the theory of  forms \cite[5.5]{W80} there is a unique self-adjoint operator $R$ in $H^2(\T)$ with $\dom R\subset \mathcal{D}_\Lambda$ and $\langle u,R\,v\rangle =\langle u,v\rangle_\Lambda$ for all $u\in \mathcal{D}_\Lambda$, $v\in\dom R$. Moreover, $R$ is determined by  $\dom R=\{v\in\mathcal{D}_\Lambda: \exists\;\widehat{v}\in H^2(\T)\textrm{ with } \langle u,v\rangle_\Lambda=\langle u,\widehat{v}\rangle\; \forall\;u\in\mathcal{D}_\Lambda\}$ and $R\,v=\widehat{v}$.  One has  $R\ge I$ and $\dom R$  dense in $\mathcal{D}_\Lambda$ with respect to $||\cdot||_\Lambda$. --- Obviously $T_\rho\subset R$. Hence $S:=R-(1-\gamma)I$ is a self-adjoint extension of $T_\omega$ with $S\ge \gamma I$. 
\\
\hspace*{6mm}
Note that   $\mathcal{D}_\omega$ in (b) is the closure  of $\dom T_\rho=\dom T_\omega$ with respect to $||\cdot||_\Lambda$, whence $\mathcal{D}_\omega \subset \mathcal{D}_\Lambda$.  Actually $\mathcal{D}_\omega =\mathcal{D}_\Lambda$ (cf. the proof of (b) in \cite[4.3]{C20}). Therefore $\dom S\subset \mathcal{D}_\omega$, which implies that $S$ is the Friedrichs extension $\tilde{T}_\omega$ of $T_\omega$.

\subsubsection*{(d)  Real square-integrable symbols} 
For square-integrable symbol $\omega$ the domain of  $T_\omega$  contains the dense subspace $\mathcal{D}_0$ of polynomials on $\T$   so that the matrix representation $(\langle e_n,T_\omega e_m\rangle)_{n,m}$ of $T_\omega|_{\mathcal{D}_0}$ is available. The latter  operator is closable since $(T_\omega|_{\mathcal{D}_0})^*\supset T^*_\omega\supset T_{\overline{\omega}}$ is densely defined. In case that $\omega$ is real and $T_\omega|_{\mathcal{D}_0}$ is essentially self-adjoint, \cite{R69} provides the spectral representation of its closure. The main results of \cite{R69} for real symbols $\omega$ regarding the spectral measure require only the integrability of $\ln^+|\omega|$, which we  recall is the necessary and sufficient condition for $T_\omega$ being densely defined.

\subsubsection*{(e)  Analytic TO}

  Let $M_\omega$ denote the operator in $H^2(\T)$ by  multiplication with a measurable $\omega$ on the domain $\dom M_\omega:=\{u\in H^2(\T):\omega u\in H^2(\T)\}$. Obviously, $M_\omega$ is closed and  is extended by $T_\omega$, i.e., $M_\omega\subset T_\omega$. --- Turn to rational $\omega=\frac{R}{S}$  with relatively prime polynomials. Then $\dom M_\omega =S_{\overline{in}}H^2(\T)$, which is dense  only if $S$ has no zeros with modulus $<1$. In this case $M_\omega=T_\omega$, which means  (identifying $H^2(\T)$ with $H^2(\mathbb{D})$ as usual) that $T_\omega$ is an analytic TO \cite{H88}, \cite[sec.\,5]{S08} (i.e., the multiplication  by $\omega$ in $H^2(\mathbb{D})$), and $\omega$ belongs to the Smirnow class $N^+$ \cite[5.2.]{S08}. If the zeros of $S$ have all modulus $>1$, then $T_\omega$ belongs to the well-studied class of bounded analytic TO (see \cite{V03} for a survey). --- For general real $\omega\in N^+$ it is shown in  \cite{H88} that $M_\omega$ has both deficiency indices finite only if $\omega$ is rational. Moreover, all pairs of deficiency indices occur, and  $M_\omega$  is self-adjoint only if $\omega$ is constant.  --- As  stated above,   $T_\omega$ and $M_\omega$ coincide  for rational symbols $\omega$ holomorphic in $\D$. More generally  $T_\omega=M_\omega$, if $\omega$ is the nontangential limit of an outer function on $\D$. Indeed, in this case $\ln|\omega|$ is integrable, whence $\ln(1+|\omega|)\le \ln2 +\ln|\omega|$ is integrable so that $T_\omega$ is densely defined (see (a)), and \cite[3.3]{CMP21} applies.

\subsubsection*{(f) Toeplitz-like operators}

 Several types of Toeplitz-like  operators have been studied in the literature. See \cite[part I]{GHJR18} for an overview. The operator $T'_\omega$ treated in \cite{GHJR18} for rational symbol $\omega$ reads $T'_\omega g=P_{H^p}f$ with $\dom T'_\omega=\{g\in H^p: \omega g=f+\rho \text{ with } f\in L^p,\,  \rho \text{ strictly proper rational with all poles on } \T\}$.  It has rather similar Fredholm properties as $T_\omega$.  Obviously $T_\omega\subset T'_\omega$ (for $p=2$). The domain of $T'_\omega$ contains  all polynomials, which allows to set up the matrix $(\langle e_n,T'_\omega e_m\rangle)_{n,m}$. It turns out to have the form of a Toeplitz matrix. Also the shift invariance $T'_\omega=T_z^{'*}\,T'_\omega \,T'_z$ holds. If $\omega$ has no poles on $\T$, i.e., $\omega$ is essentially bounded, then $\rho=0$ and $T'_\omega$ coincides with the classical Toeplitz operator $T_\omega$. Recall that $\dom T_\omega$ contains only the polynomials with factor $s$ (\ref{RSKCK})(a).


\begin{thebibliography}{9}

\bibitem{BS90} B\"ottcher A., Silbermann B.,  "Analysis of Toeplitz Operators," Springer 1990.


\bibitem{R62} Rosenblum, M., "Self-adjoint Toeplitz Operators and Associated Orthonormal Functions," Proceedings of the American Mathematical Society, vol. 13, no.4, pp. 590-595, 1962.

\bibitem{R65} Rosenblum, M.,  "A Concrete Spectral Theory for Self-adjoint  Toeplitz Operators,"  American Journal of Mathematics, vol. 87, no. 3, pp. 709-718, 1965.


\bibitem{H88} Helson, H., "Large analytic functions",  Linear Operators and Function spaces, Timi\c{s}oara, Romania,  June 6-16, 1988, pp. 209-216.  



\bibitem{S08} Sarason, D., "Unbounded Toeplitz Operators," Integr. equ. oper. theory, vol. 61,  pp. 281-298, 2008.
DOI: 10.1007/s00020-008-1588-3

\bibitem{Y18} D.R. Yafaev, D.R.  "Toeplitz versus Hankel: Semibounded Operators," Opuscula Math., vol. 38, no. 4, pp. 573-590, 2018.
 https://doi.org/10.7494/OpMath.2018.38.4.573

\bibitem{GHJR18} Groenewald, G. J.,  ter Horst, S., J. Jaftha, J., Ran, A.C.M.,  "A Toeplitz-like operator with rational symbol having poles on the unit circle, I: Fredholm properties. -- II: The Spectrum. --  III: The Adjoint",  
Operator Theory: Advances and Applications, vol. 271, pp. 239 - 268, 2018, https://doi.org/10.1007/978-3-030-04269-1--10\,\, -- Operator Theory: Advances and Applications, vol. 272, pp. 133-154, 2019. -- Integr. Equ. Oper. Theory  91:43, 2019, https://doi.org/10.1007/s00020-019-2542-2

\bibitem{C20}  Castrigiano, D.P.L., "Unbounded Wiener-Hopf Operators and Isomorphic Singular Integral Operators,"  
Complex Anal. Oper. Theory, vol.15, no. 3, 2021. DOI: https://doi.org/10.1007/s11785-021-01110-w

\bibitem{S67} Schechter, M., "Basic Theory of Fredholm Operators", Annali della Scuola Normale Superiore di Pisa, Classe di Scienze 3e s\' erie, vol.  21, no. 2, pp. 261-280, 1967.


\bibitem{CC89}  Cucker F.,  Corbalan, A.G., "An Alternate Proof of the Continuity of the Roots of a Polynomial", Amer. Math. Monthly, vol. 96, no. 4, pp. 342-345, 1989.



\bibitem{CMP21} M. C. C$\hat{a}$mara, M.T. Malheiro,  J. R. Partington: {\it Kernels of Unbounded Toeplitz Operators and Factorization of Symbols}, Results Math. \textbf{76},  no. 1, Paper No. 10 (2021)


\bibitem{C19} D.P.L. Castrigiano: {\it Spectral Representation of the Wiener-Hopf Operator for the sinc Kernel and some Related  Operators}, arXiv:2002.09235[math.FA], 2020.


\bibitem{W80}  Weidmann, J., "Linear Operators in Hilbert Spaces", Springer  1990.

\bibitem{R69}  Rovnyak, J., "On the Theory of Unbounded Toeplitz Operators", Pac. J. Math., vol. 31, no. 2, pp. 481-496, 1969.


\bibitem{V03} Vukoti\'c, D., "Analytic Toeplitz operators on the Hardy space $H^p$: a survey", Bull. Belg. Math. Soc., vol. 10, pp. 101-113, 2003.


\end{thebibliography}
\end{document}